March 8, 1996

The Characteristic Exponents of the Falling Ball Model

Nándor Simányi[1]

The Pennsylvania State University, Department of Mathematics

University Park, PA 16802

**Abstract.** We study the characteristic exponents of the Hamiltonian system of $n$ ($\geq 2$) point masses $m_1, \ldots, m_n$ freely falling in the vertical half line $\{q \mid q \geq 0\}$ under constant gravitation and colliding with each other and the solid floor $q = 0$ elastically. This model was introduced and first studied by M. Wojtkowski. Hereby we prove his conjecture: All relevant characteristic (Lyapunov) exponents of the above dynamical system are nonzero, provided that $m_1 \geq \cdots \geq m_n$ (i. e. the masses do not increase as we go up) and $m_1 \neq m_2$.

1. Introduction

In his paper [W(1990)-I] M. Wojtkowski introduced the following Hamiltonian dynamical system with discontinuities: There is a vertical half line $\{q \mid q \geq 0\}$ given and $n$ ($\geq 2$) point particles with masses $m_1 \geq m_2 \geq \cdots \geq m_n > 0$ and positions $0 \leq q_1 \leq q_2 \leq \cdots \leq q_n$ are moving on this half line so that they are subjected to a constant gravitational acceleration $a = -1$ (they fall down), they collide elastically

---

[1]On leave from the Mathematical Institute of the Hungarian Academy of Sciences, H-1364, Budapest, P. O. B. 127, Hungary

with each other, and the first (lowest) particle also collides elastically with the hard floor $q = 0$. We fix the total energy

$$H = \sum_{i=1}^{n} \left( m_i q_i + \frac{1}{2} m_i \dot{q}_i^2 \right)$$

by taking $H = 1$. The arising Hamiltonian flow with collisions $(\mathbf{M}, \{\psi^t | t \in \mathbb{R}\}, \mu)$ ($\mu$ is the Liouville measure) is the subject of this paper.

Before formulating the result of this article, however, it is worth mentioning here three important facts:

(1) Since the phase space $\mathbf{M}$ is compact, the Liouville measure $\mu$ is finite.

(2) The phase points $x \in \mathbf{M}$ for which the orbit $\{\psi^t(x) | t \in \mathbb{R}\}$ hits at least one singularity (i. e. a multiple collision) are contained in a countable union of proper, smooth submanifolds of $\mathbf{M}$ and, therefore, such points form a set of $\mu$ measure zero.

(3) For $\mu$-almost every phase point $x \in \mathbf{M}$ the collision moments of the orbit $\{\psi^t(x) | t \in \mathbb{R}\}$ do not have any finite accumulation point, see the Appendix.

In the paper [W(1990)-I] Wojtkowski formulated his main conjecture pertaining to the dynamical system $(\mathbf{M}, \{\psi^t | t \in \mathbb{R}\}, \mu)$:

**Wojtkowski's Conjecture.** *If $m_1 \geq m_2 \geq \cdots \geq m_n > 0$ and $m_1 \neq m_n$, then all but one characteristic (Lyapunov) exponents of the flow $(\mathbf{M}, \{\psi^t | t \in \mathbb{R}\}, \mu)$ are nonzero.*

**Remarks.**

1. The only exceptional exponent zero must correspond to the flow direction.

2. The condition of nonincreasing masses (as above) is essential for establishing the invariance of the symplectic cone field — an important condition for obtaining nonzero characteristic exponents. As Wojtkowski pointed out in Proposition 4 of [W(1990)-I], if $n = 2$ and $m_1 < m_2$, then there exists a linearly stable periodic orbit, thus dimming the chances of proving ergodicity.

In pursuing the goal of proving this conjecture, Wojtkowski obtained the following results in [W(1990)-I]:

**Proposition 1.** *For every $\epsilon > 0$ there is a $\delta > 0$ such that if $m_1 > m_2 > \cdots > m_n > 0$ and $\dfrac{m_1 - m_n}{m_1} < \delta$, then $(\mathbf{M}, \{\psi^t | t \in \mathbb{R}\}, \mu)$ has exactly one zero characteristic exponent except possibly on a set of $\mu$ measure $< \epsilon$.*

**Proposition 2.** *If there are exactly $l$ groups of particles with equal masses, $l \geq 2$, containing $k_1, \ldots, k_l$ particles respectively, the greatest common divisor of $k_1, \ldots, k_l$ is one and $m_1 \geq m_2 \geq \cdots \geq m_n > 0$, then $\{\psi^t\}$ has exactly one zero characteristic exponent on a set of positive Liouville measure.*

**Proposition 3.** *If $n = 3$ and $m_1 > m_2 > m_3$, then $\{\psi^t\}$ has exactly one zero characteristic exponent $\mu$-almost everywhere.*

In the subsequent article [W(1990)-II] Wojtkowski replaced the linear potential $U(q) = q$ of constant gravitation by a varying gravitational force with potential $U(q)$ for which $U'(q) > 0$ and $U''(q) < 0$. (The usual gravitational potential $U(q) = \dfrac{-1}{q + c_0}$ belongs to this category!) He proved there that in the falling ball system with such a potential $U(q)$ all relevant characteristic exponents are nonzero almost everywhere.

The result of this paper is a slightly weakened version of Wojtkowski's conjecture:

**Theorem.** *If $m_1 > m_2 \geq m_3 \geq \cdots \geq m_n > 0$, then $\mu$-almost everywhere all but one characteristic exponents of the flow $(\mathbf{M}, \{\psi^t | t \in \mathbb{R}\}, \mu)$ are nonzero.*

We are closing this brief introduction by mentioning that in his work [Ch(1993)] N. I. Chernov significantly relaxed a condition of the Liverani–Wojtkowski local ergodicity theorem for symplectomorphisms, [L-W(1995)]. (This theorem is a generalization of the celebrated Theorem on Local Ergodicity for semi-dispersing billiards by Chernov and Sinai, [S-Ch(1987)].) The ominous condition is the "proper alignment" of the singularity manifolds, Condition D in Section 7 of [L-W(1995)]. This condition is easily seen to be violated by the system of falling balls (see Section 14.C of [L-W(1995)]), but the relaxed condition 5' of Chernov's paper [Ch(1993)] is very likely to hold for the system of falling balls. Thus, thanks to Chernov's im-provement in the local ergodicity theorem, the ergodicity of the falling ball system essentially got within our reach!

For a more detailed introduction to this subject, and for a thoroughly assembled collection of references and historical remarks, the reader is kindly referred to the introduction of the paper [W(1990)-I].

**Remark.** It is easy to see that the studying of the falling point particles on the vertical half line $\{q \mid q \geq 0\}$ is not a restriction of generality as compared to the systems of 1-D balls (hard rods) of length $2r$. Namely, the simple change in the kinetic data $q_i \mapsto q_i - (2i-1)r$, $v_i \mapsto v_i$, $H_0 \mapsto H_0 - r\sum_{i=1}^n (2i-1)m_i$ (the change of the fixed level of energy) establishes an isomorphism between the hard rod system and the point particle model.

## 2. Prerequisites

The upcoming brief survey of our dynamical system and the related technicalities will narrowly follow the approach of sections 1–3 of [W(1990)-I]. A thorough description of the falling ball model and detailed references can be found in that article.

We consider the following Hamiltonian system with discontinuities: Given the vertical half line $\{q \mid q \geq 0\}$, $n$ ($\geq 2$) point particles (one dimensional "balls" labelled by $1, 2, \ldots, n$) with positions $0 \leq q_1 \leq q_2 \leq \cdots \leq q_n$ are moving in that half line so that they fall down under a constant gravitational acceleration $a = -1$, they collide with each other elastically whenever they hit each other, and the first particle bounces back from the hard floor $q = 0$ elastically when it hits the floor. Denote by $v_i = \dot{q}_i$ the velocity, by $p_i = m_i v_i$ the momentum, and by $h_i = m_i q_i + \frac{1}{2} m_i v_i^2 = m_i q_i + \frac{p_i^2}{2m_i}$ the energy of the i-th ball. (The quantity $m_i q_i$ is the potential energy arising from the gravitation.) The extended phase space (without fixing the energy)

4of this mechanical system is then

(2.1) $$\mathbf{N} = \{(q,p) \in \mathbb{R}^n \times \mathbb{R}^n \mid 0 \leq q_1 \leq \cdots \leq q_n\}.$$

The manifold (with boundary) $\mathbf{N}$ carries the usual symplectic form

(2.2) $$\omega = \sum_{i=1}^{n} dq_i \wedge dp_i = \sum_{i=1}^{n} dh_i \wedge dv_i.$$

The Hamiltonian function is $H(q,p) = \sum_{i=1}^{n} h_i$, and the arising Hamiltonian flow $\{\psi^t \mid t \in \mathbb{R}\}$ is determined by the usual formalism

(2.3) $$\begin{cases} \dot{q}_i = \dfrac{p_i}{m_i} = \dfrac{\partial H}{\partial p_i} \\ \dot{p}_i = -m_i = -\dfrac{\partial H}{\partial q_i} \end{cases}$$

$(i = 1, \ldots, n)$ between collisions.

A collision of type $(i, i+1)$ $(i = 1, \ldots, n-1)$ occurs when $q_i = q_{i+1}$. Then the velocities and the momenta of the colliding particles get transformed according to the law of elastic collisions:

(2.4) $$\begin{cases} v_i^+ = \gamma_i v_i^- + (1 - \gamma_i) v_{i+1}^- \\ v_{i+1}^+ = (1 + \gamma_i) v_i^- - \gamma_i v_{i+1}^- \\ p_i^+ = \gamma_i p_i^- + (1 + \gamma_i) p_{i+1}^- \\ p_{i+1}^+ = (1 - \gamma_i) p_i^- - \gamma_i p_{i+1}^-, \end{cases}$$

where $\gamma_i = \dfrac{m_i - m_{i+1}}{m_i + m_{i+1}}$, see also (3) of [W(1990)-I]. (Here the superscript $-$ $(+)$ refers to the kinetic data measured right before (after) the considered collision.)

At a floor collision $q_1 = 0$ we obviously have

(2.5) $$\begin{cases} v_1^+ = -v_1^- \\ p_1^+ = -p_1^-. \end{cases}$$

We introduce the following notations for the several components of the boundary $\partial \mathbf{N}$ of $\mathbf{N}$:

(2.6) $$\begin{cases} \partial \mathbf{N}_i^\pm = \{(q,p) \in \mathbf{N} \mid q_i = q_{i+1} \ \& \ \pm(v_i - v_{i+1}) < 0\} \\ \partial \mathbf{N}^\pm = \cup_{i=0}^{n-1} \partial \mathbf{N}_i^\pm \end{cases}$$

for $i = 0, \ldots, n-1$, where, by convention, $q_0 = v_0 = 0$. Then the collision map $\Phi_i : \partial \mathbf{N}_i^- \to \partial \mathbf{N}_i^+$ ($i = 0, \ldots, n-1$), determined by (2.4)–(2.5) and by the condition that the positions do not change at collisions, turns out to be a symplectic diffeomorphism of $\partial \mathbf{N}_i^-$ onto $\partial \mathbf{N}_i^+$, see Section 2 of [W(1990)-I]. (Here the $(2n-1)$-dimensional manifold $\partial \mathbf{N}$ naturally inherits the pseudo-symplectic structure $\omega \lceil \partial \mathbf{N}$, i. e. the restriction of $\omega$ to $\partial \mathbf{N}$.)

Now we fix the total energy by taking $H = H_0 = 1$, and consider the restriction of the Hamiltonian flow with collisions $\{\psi^t \mid t \in \mathbb{R}\}$ to the energy hypersurface

$$\mathbf{M} = \{(q,p) \in \mathbf{N} \mid H(q,p) = 1\}.$$

The corresponding boundary components of $\mathbf{M}$ are denoted by $\partial \mathbf{M}_i^\pm$ and $\partial \mathbf{M}^\pm$.

The (invariant) Liouville measure $\nu$ in the extended phase space $\mathbf{N}$ is defined via the volume element $\mathcal{V} = \omega \wedge \cdots \wedge \omega$ (the $n$-th exterior power of $\omega$). The corresponding conditional Liouville measure $\mu$ on $\mathbf{M}$ can then be obtained as the contraction $\iota_F(\mathcal{V})$ of the $2n$-dimensional volume element $\mathcal{V}$ by a vector field $\{F(x) \mid x \in \mathbf{M}\}$ for which $D_F(H) = 1$. Under our assumptions the phase space $\mathbf{M}$ is compact and, therefore, the measure $\mu$ is a finite, $\{\psi^t\}$-invariant Borel measure on $\mathbf{M}$. The subject of this paper is the Hamiltonian flow with collisions $(\mathbf{M}, \{\psi^t\}, \mu)$.

*Factorization with respect to the flow direction*

We will frequently use another coordinate system $(\delta h, \delta v)$ in the tangent space $\mathcal{T}_x \mathbf{N}$:

(2.7) $$\begin{cases} \delta h_i = m_i \delta q_i + v_i \delta p_i \\ \delta v_i = \dfrac{1}{m_i} \delta p_i. \end{cases}$$

For every interior point $x = (q,p)$ of $\mathbf{M}$ we define the codimension-one linear subspace $\mathcal{T}_x$ of the tangent space $\mathcal{T}_x \mathbf{M}$ of $\mathbf{M}$ at $x$ as follows:

(2.8) $$\mathcal{T}_x = \left\{(\delta h, \delta v) \in \mathcal{T}_x \mathbf{M} \,\bigg|\, \sum_{i=1}^n \delta v_i = 0\right\}.$$

It is clear that the subspace $\mathcal{T}_x$ of $\mathcal{T}_x\mathbf{M}$ is transversal to the velocity vector

$$V(x) = \frac{d}{dt}\psi^t(x)\bigg|_{t=0} = (0,\ldots,0;-1,\ldots,-1)$$

(in $(\delta h, \delta v)$ coordinates) of the flow $\{\psi^t | t \in \mathbb{R}\}$.

Since the $\omega$-orthocomplement of the tangent vector $V(x)$ ($x \in \text{int}\mathbf{M}$) in $\mathcal{T}_x\mathbf{N}$ is precisely the space $\mathcal{T}_x\mathbf{M}$, we infer that the 2-form $\omega\lceil\mathcal{T}_x\mathbf{M}$ naturally descends to the factor space

$$\mathcal{T}_x\mathbf{M}/\{\lambda V(x)|\, \lambda \in \mathbb{R}\} = \mathcal{T}_x\mathbf{M}/V_x,$$

and it is a nondegenerate 2-form on that space. Moreover, the composition of the inclusion $\mathcal{T}_x \hookrightarrow \mathcal{T}_x\mathbf{M}$ and the projection $\mathcal{T}_x\mathbf{M} \to \mathcal{T}_x\mathbf{M}/V_x$ provides a natural identification $\iota : \mathcal{T}_x \to \mathcal{T}_x\mathbf{M}/V_x$. The linearization (derivative) $D\psi^t(x) : \mathcal{T}_x\mathbf{M} \to \mathcal{T}_{\psi^t(x)}\mathbf{M}$ maps the line $V_x$ onto $V_{\psi^t(x)}$ and, therefore, it descends to a mapping

$$D\psi^t(x) : \mathcal{T}_x\mathbf{M}/V_x \longrightarrow \mathcal{T}_{\psi^t(x)}\mathbf{M}/V_{\psi^t(x)}.$$

Thus, by using the abovementioned identification $\iota$, we can (and will) think of $D\psi^t(x)$ as a mapping $D\psi^t(x) : \mathcal{T}_x \to \mathcal{T}_{\psi^t(x)}$. The facts that the collision maps $\Phi_i : \partial\mathbf{M}_i^- \to \partial\mathbf{M}_i^+$ are symplectomorphisms and $\{\psi^t | t \in \mathbb{R}\}$ is a Hamiltonian flow between collisions together imply that the linearization $D\psi^t(x) : \mathcal{T}_x \to \mathcal{T}_{\psi^t(x)}$ preserves the nondegenerate 2-form $\omega$.

*Characteristic (Lyapunov) exponents*

In the manifold $\mathbf{N}$ we introduce the Riemannian metric compatible with the 2-form $\omega = \sum_{i=1}^n dh_i \wedge dv_i$:

(2.9) $$\|(\delta h, \delta v)\|^2 = \sum_{i=1}^n \left[(\delta h_i)^2 + (\delta v_i)^2\right],$$

$(\delta h, \delta v) \in \mathcal{T}_x\mathbf{N}$. It is worth noting here that there could be several other natural possibilities of defining a Riemannian metric in $\mathbf{N}$, for example by using $\|\delta q\|^2$ and $\|\delta p\|^2$ instead of $\|\delta h\|^2$ and $\|\delta v\|^2$. However, it is an easy consequence of (2.7) that the ratio between these metrics would remain between two positive constants so that those constants would only depend on the masses of the particles and, being so, such a change in the Riemannian metric would not affect the characteristic exponents. It is obvious that the space $\mathcal{T}_x$ ($x \in \text{int}\mathbf{M}$) is just the orthogonal complement (with respect to the Riemannian metric of (2.9)) of the line $V_x$ in the space $\mathcal{T}_x\mathbf{M}$.

Thanks to the arising inner product in the space $\mathcal{T}_x$, we can now speak about the adjoint linear map

$$\left(D\psi^t(x)\right)^* : \mathcal{T}_{\psi^t(x)} \longrightarrow \mathcal{T}_x$$

of the linearization $D\psi^t(x) : \mathcal{T}_x \to \mathcal{T}_{\psi^t(x)}$. Both mappings are symplectic (and, therefore, invertible) with respect to the restriction of the 2-form $\omega$.

We will be formulating Oseledets' Multiplicative Ergodic Theorem for a suitable Poincaré section $\Psi$ (a discrete time version) of the flow $(\mathbf{M}, \{\psi^t\}, \mu)$. Namely, we take $\Psi : \partial\mathbf{M}^+ \to \partial\mathbf{M}^+$ the first return map to the boundary $\partial\mathbf{M}^+$ with respect to the flow $\{\psi^t\}$. The natural invariant measure $\mu_0$ of $\Psi$ on $\partial\mathbf{M}^+$ can be obtained by taking $\mu_0(A) := \lim_{\epsilon \searrow 0} \frac{1}{\epsilon}\mu(A_\epsilon)$, $A \subset \partial\mathbf{M}^+$, where

$$A_\epsilon = \left\{x \in \mathbf{M} | \ \exists y \in A \ \exists t, \ 0 \leq t \leq \epsilon \text{ such that } x = \psi^t(y)\right\}.$$

It is a classically known result that the mapping $(\partial\mathbf{M}^+, \Psi, \mu_0)$ has (almost everywhere) nonzero Lyapunov exponents if and only if the flow $(\mathbf{M}, \{\psi^t\}, \mu)$ has the similar property. The reason is that the Lyapunov exponents of the map $\Psi$ are just the corresponding characteristic exponents of the flow multiplied by the average (with respect to $\mu_0$) return time $t(x)$ ($x \in \partial\mathbf{M}^+$) to the boundary $\partial\mathbf{M}^+$, see, for instance, Lemma 2.2 of [W(1985)]. But the ergodic average $t(x)$ is a $\mu_0$-almost everywhere positive function.

We introduce the shorthand $\mathcal{L}_x = \mathcal{T}_x\partial\mathbf{M}^+$ for the points $x \in \partial\mathbf{M}^+$. The tangent space $\mathcal{L}_x$ is a Euclidean space with the restriction of the inner product (2.9). The linearization $D\Psi(x) : \mathcal{L}_x \to \mathcal{L}_{\Psi(x)}$ and its adjoint $(D\Psi(x))^* : \mathcal{L}_{\Psi(x)} \to \mathcal{L}_x$ are now symplectic isomorphisms between $\mathcal{L}_x$ and $\mathcal{L}_{\Psi(x)}$. The integrability condition

$$\int_{\partial \mathbf{M}^+} \|D\Psi(x)\| \, d\mu_0(x) < \infty$$

of the Oseledets' theorem is an easy consequence of (10) from [W(1990)-I] and the following simple facts:

(i) $\dfrac{d}{dt}(\delta h) = \dfrac{d}{dt}(\delta v) = 0$ (between collisions),

(ii) the return time from $\partial \mathbf{M}^+$ to $\partial \mathbf{M}^+$ is bounded.

The conditions of Oseledets' celebrated Multiplicative Ergodic Theorem, [O(1968)], [R(1979)], are fulfilled by the invertible measurable cocycle $(\partial \mathbf{M}^+, \Psi, D\Psi, \mu_0)$:

**Multiplicative Ergodic Theorem.** *For $\mu$-almost every phase point $x \in \partial \mathbf{M}^+$ the following limit exists:*

$$\lim_{n \to \infty} \left[ (D\Psi^n(x))^* \circ D\Psi^n(x) \right]^{\frac{1}{2n}} := \Lambda(x).$$

*The invertible limit mapping $\Lambda(x) : \mathcal{L}_x \to \mathcal{L}_x$ is a symmetric, positive, symplectic linear transformation. The logarithms of the eigenvalues of $\Lambda(x)$*

$$-\lambda_{n-1}(x) \leq -\lambda_{n-2}(x) \leq \cdots \leq -\lambda_1(x) \leq 0 \leq \lambda_1(x) \leq \cdots \leq \lambda_{n-1}(x)$$

*are called the characteristic (Lyapunov) exponents of the invertible measurable cocycle $(\partial \mathbf{M}^+, \Psi, D\Psi, \mu_0)$.*

*The invariant cone field*

The symplectic linear space $\mathcal{T}_x \mathbf{N}$ is the direct sum $V_1 \oplus V_2$ of the Lagrangian subspaces

$$\begin{cases} V_1 = \{(\delta h, \delta v) \in \mathcal{T}_x \mathbf{N} | \, \delta v = 0\} \\ V_2 = \{(\delta h, \delta v) \in \mathcal{T}_x \mathbf{N} | \, \delta h = 0\} . \end{cases}$$

In Section 4 of [L-W(1995)] Liverani and Wojtkowski introduced the nondegenerate quadratic form $Q$ in $\mathcal{T}_x \mathbf{N}$ associated with the decomposition $\mathcal{T}_x \mathbf{N} = V_1 \oplus V_2$ as follows:

(2.10) $$Q((\delta h, \delta v)) = \langle \delta h; \delta v \rangle = \sum_{i=1}^{n} \delta h_i \delta v_i.$$

The corresponding positive cone (sector) $C_x \subset \mathcal{T}_x \mathbf{N}$ between the Lagrangian subspaces $V_1$ and $V_2$ is then defined as follows:

$$(2.11) \qquad C_x = \{(\delta h, \delta v) \in \mathcal{T}_x \mathbf{N} | \, Q\left((\delta h, \delta v)\right) \geq 0\}.$$

We will use the restriction of the quadratic form $Q$ to the space $\mathcal{T}_x$ and the intersection $C_x \cap \mathcal{T}_x$, also (a bit sloppily) denoted by $Q$ and $C_x$. It is worth noting here that the $Q$-orthocomplement of the line $V_x$ in $\mathcal{T}_x \mathbf{N}$ is precisely the space $\mathcal{T}_x \mathbf{M}$ and, therefore, the form $Q$ descends to a nondegenerate quadratic form (also denoted by $Q$) on the factor space $\mathcal{T}_x \mathbf{M}/V_x \cong \mathcal{T}_x$.

Wojtkowski shows in Section 4 of [W(1990)-I] that if $m_1 \geq m_2 \geq \cdots \geq m_n$, then for every $t > 0$ the linearized mapping $D\psi^t(x) : \mathcal{T}_x \to \mathcal{T}_{\psi^t(x)}$ is $Q$-monotonic, i. e. for every tangent vector $y \in \mathcal{T}_x$ one has $Q\left[(D\psi^t(x))(y)\right] \geq Q(y)$ or, equivalently, the cone field $C$ is invariant:

$$(2.12) \qquad \left(D\psi^t(x)\right)(C_x) \subset C_{\psi^t(x)}$$

for every $t > 0$, $x \in \text{int}\mathbf{M}$, $\psi^t(x) \in \text{int}\mathbf{M}$.

**Definition 2.13.** *Let $\{\psi^t(x) | \, t \geq 0\}$ be a nonsingular, positive orbit, $x \in \text{int}\mathbf{M}$. We say that the cone field $C$ is **eventually strictly invariant** along the orbit $\{\psi^t(x) | \, t \geq 0\}$ iff there exits a number $t_0 > 0$ such that*

$$\left(D\psi^{t_0}(x)\right)(C_x) \subset \text{int}\left(C_{\psi^{t_0}(x)}\right),$$

*or, equivalently, $Q\left[(D\psi^{t_0}(x))(y)\right] > 0$ for every $0 \neq y \in C_x$.*

A major result (Theorem 5.1) of [W(1985)] is the following one:

**Theorem on Nonzero Characteristic Exponents.** *If the cone field $C$ is eventually strictly invariant along $\{\psi^t(x) | \, t \geq 0\}$ for $\mu$-almost every $x \in \mathbf{M}$, then all characteristic exponents $\lambda_i(x)$ of the cocycle $(\partial \mathbf{M}^+, \Psi, D\Psi, \mu_0)$ are different from zero for $\mu_0$-almost every $x \in \partial \mathbf{M}^+$.*

In Section 3 we will just check the conditions of this theorem for the falling ball system introduced before.

## 3. The Strict Invariance of the Cone Field

### (Proof of the Theorem)

In this section we will be studying the strict cone invariance along a non-singular trajectory

$$\{\psi^t(x)|\, t \in \mathbb{R}\} = \{(q_1(t), \ldots, q_n(t); v_1(t), \ldots, v_n(t))\,|\, t \in \mathbb{R}\}$$

of the Hamiltonian flow with collisions $\{\psi^t|\, t \in \mathbb{R}\}$ introduced in Section 1. We will always assume that $m_1 > m_2 \geq m_3 \geq \cdots \geq m_n > 0$ and $t = 0$ is not a moment of collision. Through references it was shown in the previous section that the quantities $(\delta h; \delta v) = (\delta h_1, \ldots, \delta h_n; \delta v_1, \ldots, \delta v_n)$ (where we always assume $\sum_{i=1}^n \delta h_i = \sum_{i=1}^n \delta v_i = 0$) serve as suitable symplectic coordinates in the codimension-one subspace $\mathcal{T}_x$ of the tangent space $\mathcal{T}_x\mathbf{M}$ of $\mathbf{M}$ at the phase point $x = (q_1(0), \ldots, q_n(0); v_1(0), \ldots, v_n(0))$. Recall that the linear space $\mathcal{T}_x$ is transversal to the flow direction and the restriction of the canonical symplectic form

$$(3.1) \qquad \omega = \sum_{i=1}^n dh_i \wedge dv_i = \sum_{i=1}^n dq_i \wedge dp_i$$

of $\mathbf{M}$ is nondegenerate on $\mathcal{T}_x$. We also recall from the previous section that the individual energy of the i-th particle is $h_i = m_i q_i + \frac{1}{2} m_i v_i^2 = m_i q_i + \frac{p_i^2}{2m_i}$ and, therefore, $\delta h_i = m_i \delta q_i + m_i v_i \delta v_i = m_i \delta q_i + v_i \delta p_i$. In his article [W(1990)-I] Wojtkowski introduced the nondegenerate quadratic form $Q\left((\delta h; \delta v)\right) = \langle \delta h; \delta v \rangle$ (the usual Euclidean inner product of $\delta h$ and $\delta v$) and the corresponding positive cone field

$$(3.2) \qquad C_x = \{(\delta h; \delta v) \in \mathcal{T}_x |\, Q\left((\delta h; \delta v)\right) \geq 0\},$$

see also the introduction of the present paper. He proved in the mentioned article that the cone field $C_x$ is invariant under the linearization of the flow $\{\psi^t|\, t \in \mathbb{R}\}$, i. e. $D\psi^t(C_x) \subset C_{\psi^t(x)}$ for every $t \geq 0$. (For a detailed discussion of the algebra and geometry of an invariant cone field, see §4-6 of [L-W(1995)].) Furthermore, Liverani and Wojtkowski also proved in §4-6 of [L-W(1995)] that the *eventually strict invariance* of the cone field (i. e. $D\psi^t(C_x) \subset \text{int}\left(C_{\psi^t(x)}\right)$ for some $t = t(x) > 0$ and for $\mu$-almost every $x \in \mathbf{M}$) implies that all but one characteristic exponents of the flow $\{\psi^t | t \in \mathbb{R}\}$ are nonzero! The only exceptional zero exponent corresponds to the flow direction.

It follows from Wojtkowski's arguments between the Theorem and Proposition 1 of Section 5 of [W(1990)-I] that in order to check the eventually strict invariance of the cone field along the studied nonsingular trajectory $\{\psi^t(x) | t \in \mathbb{R}\}$ it is enough to prove that

(A) for every vector $0 \neq (0; \delta v) \in \mathcal{T}_x$ there exists a $t > 0$ such that

$$Q\left[D\psi^t\left((0; \delta v)\right)\right] > 0$$

and

(B) for every vector $0 \neq (\delta h; 0) \in \mathcal{T}_x$ there exists a $t > 0$ such that

$$Q\left[D\psi^t\left((\delta h; 0)\right)\right] > 0.$$

Moreover, Wojtkowski's mentioned arguments from Section 5 of [W(1990)-I] actually contain the proof of (A), provided that $m_1 > m_2 \geq \cdots \geq m_n > 0$. Here we briefly review his ideas. The formula (13) of [W(1990)-I] says that a tangent vector of the form $(0; \delta v)$ does not get changed at all by the linearization $D\Phi_0$ of the collision map $\Phi_0$ corresponding to the floor collision. Suppose that a collision of type $(i, i+1)$ occurs at time $t_k$. Denote the corresponding collision map by $\Phi_i$ as in [W(1990)-I]. Suppose for a while that $m_i > m_{i+1}$. The equations (10) and (11) of the above article say that after pushing the tangent vector $(0; \delta v)$ through the collision $(i, i+1)$, the value of the $Q$ form on the image $D\Phi_i\left((0; \delta v)\right)$ either becomes positive, or $\delta v_i^-(t_k) = \delta v_{i+1}^-(t_k)$ and $\delta v^+(t_k) = \delta v^-(t_k)$.

On the other hand, it also follows from (10) of the mentioned paper that in the case $m_i = m_{i+1}$ the linearization $D\Phi_i$ of the collision map $\Phi_i$ simply interchanges $\delta v_i$ and $\delta v_{i+1}$: $\delta v_i^+(t_k) = \delta v_{i+1}^-(t_k)$, $\delta v_{i+1}^+(t_k) = \delta v_i^-(t_k)$, $\delta v_j^+(t_k) = \delta v_j^-(t_k)$ for $j \ne i, i+1$, and $\delta h^+ = \delta h^- = 0$. Therefore, it is quite reasonable to re-label the particles dynamically at every such collision by simply interchanging the labels $i$ and $i+1$. This is equivalent to allowing the particles with equal masses to freely penetrate through each other precisely the same way as Wojtkowski did in [W(1990)-I]. Then, as long as the $Q$ form remains zero on the images of $(0; \delta v)$, the images of $(0; \delta v)$ under the linearization of the flow $\{\psi^t | t \in \mathbb{R}\}$ remain the same, and $\delta v_i = \delta v_j$ if the particles $i$ and $j$ collide on the considered trajectory segment. Since each particle $i$ with $m_i < m_1$ must eventually bounce back from a heavier particle, and 1 is the sole heaviest particle by our assumption $m_1 > m_2$, we obtain that every $\delta v_i$ must be the same. By our convention $\sum_{i=1}^n \delta v_i = 0$ (we are always dealing with vectors from $\mathcal{T}_x$), however, we infer that $\delta v = 0$, provided that all future images of the considered tangent vector $(0; \delta v)$ have zero $Q$ form.

Thus, in order to prove the Theorem, it is enough to show that (B) holds true for $\mu$-almost every $x \in \mathbf{M}$. This is what we are going to do.

We begin with the definition of the "neutral space" $\mathcal{N}_x$ of the nonsingular phase point $x$. (To be more accurate, $\mathcal{N}_x$ is going to be the neutral space of the positive orbit $\{\psi^t(x) | t \ge 0\}$.) The linear subspace $\mathcal{N}_x$ of $\mathcal{T}_x$ will be the precise analogue of the neutral space $\mathcal{N}_0\left(S^{[0,\infty]}x\right)$ of a positive orbit in a semi-dispersing billiard, originally and essentially introduced by Chernov and Sinai in [S-Ch(1987)], and later heavily used by Krámli, Szász and myself in the several proofs of ergodicity for hard ball systems, [K-S-Sz(1991)], [K-S-Sz(1992)], [Sim(1992)-I-II], [S-Sz(1995)], [S-Sz(1996)].

**Definition 3.3.**

$$\mathcal{N}_x := \left\{(\delta h; 0) \in \mathcal{T}_x \,|\, Q\left[D\psi^t\left((\delta h; 0)\right)\right] = 0 \quad \forall t \ge 0\right\}.$$

It is easy to convince ourselves that, indeed, $\mathcal{N}_x$ is a linear subspace of $\mathcal{T}_x$. The main result of this section (immediately proving the theorem) is

**Main Lemma 3.4.** *For $\mu$-almost every (nonsingular) phase point $x \in \mathbf{M}$ we have $\mathcal{N}_x = \{0\}$.*

**Proof.** The proof will be based on a few lemmas. Denote by $0 < t_1 < t_2 < \ldots$, $t_k \nearrow \infty$ the sequence of all collision moments on the positive orbit $\{\psi^t(x) | t \geq 0\}$. (These collision moments do not have a finite accumulation point, see the Appendix.) It is obvious that the tangent vector $D\psi^t((\delta h; 0)) := (\delta h(t); \delta v(t))$ does not change between collisions (after the natural identification of the tangent spaces of $\mathbf{M}$ at different points). Between the Theorem and Proposition 1 of Section 5 of [W(1990)-I] Wojtkowski proved that for any neutral vector $(\delta h; 0) \in \mathcal{N}_x$ we necessarily have that

(1) $\delta v(t) = 0$, i. e. $D\psi^t((\delta h; 0)) = (\delta h(t); 0) \quad \forall t \geq 0$;

(2) $\delta h^-(t_k) = \delta h^+(t_k)$ and $\delta h_1(t_k) = 0$ if the first particle bounces back from the floor at time $t_k$, i. e. $q_1(t_k) = 0$;

(3) $\delta h^+(t_k) = R_i^* \delta h^-(t_k)$ if $t_k$ is a moment of an $(i, i+1)$ collision $(1 \leq i \leq n-1)$, where $R_i$ is the following $n \times n$ matrix:

$$(3.5) \quad R_i = \begin{bmatrix} 1 & 0 & \ldots & & \ldots & 0 & 0 \\ 0 & \ddots & \vdots & & \vdots & 0 & 0 \\ 0 & \ldots & \gamma_i & 1-\gamma_i & \ldots & 0 \\ 0 & \ldots & 1+\gamma_i & -\gamma_i & \ldots & 0 \\ 0 & 0 & \vdots & & \vdots & \ddots & 0 \\ 0 & 0 & \ldots & & \ldots & 0 & 1 \end{bmatrix}.$$

Here $\gamma_i = \dfrac{m_i - m_{i+1}}{m_i + m_{i+1}}$ and the four entries containing $\gamma_i$ fill up the intersections of the i-th and j-th rows and columns. (For property (3), see also (10) of [W(1990)-I].)

**Lemma 3.6.** *For every vector $(\delta h; 0) \in \mathcal{N}_x$ we have $\sum_{i=1}^n \delta h_i v_i(0) = 0$ and, hence, $\sum_{i=1}^n \delta h_i(t) v_i(t) = 0$ for all $t \geq 0$.*

**Proof of Lemma 3.6.** Consider the quantity

$$w(t) := \sum_{i=1}^n q_i(t) \delta h_i(t)$$

well defined for all $t \geq 0$. The obvious relation $\dfrac{d}{dt}\delta h_i(t) = 0$ implies that *between collisions*

$$\text{(3.7)} \qquad \frac{d}{dt}w(t) = \sum_{i=1}^{n} v_i(t)\delta h_i(t) = \langle v(t); \delta h(t)\rangle$$

and

$$\text{(3.8)} \qquad \frac{d^2}{dt^2}w(t) = -\sum_{i=1}^{n} \delta h_i(t) = 0.$$

(Here we took advantage of the fact that $\dot v_i = -1$ is the gravitational acceleration.) Thus $w(t)$ is a linear function of $t$ between collisions. It is a straightforward consequence of (2) that the function $w$ is continuous and it even does not change its slope at a collision with the floor. If, on the other hand, a collision of type $(i, i+1)$ takes place at time $t_k$ ($1 \leq i \leq n-1$), then we have that $q_i(t_k) = q_{i+1}(t_k)$ and the compound velocity vector $v(t)$ gets transformed by the matrix $R_i$ at time $t_k$: $v^+(t_k) = R_i v^-(t_k)$, see (9) of [W(1990)-I]. These facts and property (3) imply that $w^+(t_k) = w^-(t_k)$ and

$$\frac{d}{dt}w(t_k + 0) = \langle \delta h(t_k+0); v(t_k+0)\rangle = \langle R_i^*\delta h(t_k-0); R_i v(t_k-0)\rangle =$$

$$\langle \delta h(t_k-0); R_i^2 v(t_k-0)\rangle = \frac{d}{dt}w(t_k-0).$$

Here we used the obvious relation $R_i^2 = 1$. We have seen therefore that $w(t)$ is an (inhomogeneous) linear function of $t \geq 0$.

**Sublemma 3.9.** *The function $w(t)$ is bounded ($t \geq 0$).*

**Proof.** It is enough to prove that the quantity $\sum_{i=1}^{n} \dfrac{1}{m_i}[\delta h_i(t)]^2$ is a constant function of $t$, since the positions $q_i(t)$ are obviously bounded. However, the quantity

$$\sum_{i=1}^{n} \frac{1}{m_i}[\delta h_i(t)]^2 = \sum_{i=1}^{n} \delta h_i(t)\delta q_i(t) = \langle \delta h(t); \delta q(t)\rangle$$

is constant between collisions and, by (2), it does not change its value at a floor collision. Furthermore, if a collision of type $(i, i+1)$ takes place at $t_k$, then the vector $\delta q(t) = D^{-1}\delta h(t)$ – where $D = \text{diag}(m_1, \ldots, m_n)$ is the diagonal matrix with the masses as entries – gets transformed by the matrix $DR_i^* D^{-1} = R_i$, i. e. $\delta q(t_k + 0) = R_i \delta q(t_k - 0)$. Therefore, according to (3), we see that

$$\langle \delta h(t_k + 0); \delta q(t_k + 0) \rangle = \langle R_i^* \delta h(t_k - 0); R_i \delta q(t_k - 0) \rangle = \langle \delta h(t_k - 0); \delta q(t_k - 0) \rangle.$$

(The above arguments are precisely the arithmetic background of the conservation of the kinetic energy at an $(i, i+1)$ collision.)

Hence the sublemma follows. $\square$

The assertion of the sublemma, together with (3.7), now proves Lemma 3.6. $\square$

*Finishing the Proof of Main Lemma 3.4*

Set

(3.10) $$X_d = \{x \in \mathbf{M} | \, x \text{ is nonsingular \& } \dim \mathcal{N}_x = d\}$$

for $d = 0, 1, \ldots, n - 1$. We want to show that for every $d > 0$ the set $X_d$ has $\mu$ measure zero. Fix a number $d > 0$ and an arbitrary phase point $x_0 \in X_d$. We will prove that $x_0$ has a suitably small, open neighborhood $U_0$ in $\mathbf{M}$ with $\mu(X_d \cap U_0) = 0$.

First choose a number $\tau > 0$ with the following properties:

(i) $\tau$ is not a moment of collision in the positive trajecory of $x_0$;

(ii) $\mathcal{N}_{x_0} = \mathcal{N}(\{\psi^t(x_0) | \, 0 \leq t \leq \tau\})$, i. e. the relation ($\forall t, 0 \leq t \leq \tau$) $Q[D\psi^t((\delta h; 0))] = 0$ implies that ($\forall t \geq 0$) $Q[D\psi^t((\delta h; 0))] = 0$.

It is a very important consequence of properties (1)-(3) above that *the neutral space of a finite orbit segment is completely determined by the symbolic collision sequence of that segment, i. e. by the types of collisions!* Therefore, one can surely find a small, open neighborhood $U_0$ of $x_0$ in $\mathbf{M}$ such that for every element $y \in U_0$ we have that

(i)' the symbolic collision sequence of $\{\psi^t(y)|\, 0\le t\le \tau\}$ is the same as of

$$\{\psi^t(x_0)|\, 0\le t\le \tau\}$$

and $\tau$ is not a collision moment of these orbit segments;

(ii)' $\mathcal{N}\left(\{\psi^t(y)|\, 0\le t\le \tau\}\right) = \mathcal{N}_{x_0}$.

Then we conclude that
$$\begin{cases} \dim \mathcal{N}_{x_0} = d \quad (\ge 1), \\ (\forall y \in U_0)\ \mathcal{N}_y \subset \mathcal{N}_{x_0}. \end{cases}$$

Therefore, we get that $\mathcal{N}_y = \mathcal{N}_{x_0}$ for every $y \in U_0 \cap X_d$. According to Lemma 3.6, however, the compound velocity $v_y$ of such a point $y = (q_y, v_y) \in U_0 \cap X_d$ is necessarily perpendicular to the $\delta h$ part of every neutral vector $(\delta h; 0) \in \mathcal{N}_{x_0}$. Since $\mathcal{N}_{x_0} \ne \{0\}$ and the velocities $v_y$ of all points $y \in U_0$ fill out an open subset of $\mathbb{R}^n$, we see that, indeed, $\mu(U_0 \cap X_d) = 0$.

This finishes the proof of Main Lemma 3.4 and, therefore, the proof of the Theorem, as well. $\square$

## APPENDIX

### Degenerate Orbits

We begin with a simple observation. Assume that we are given a phase point $x = (q, p) \in \mathbf{M}$ with the property that there are $k$ particles ($1 \le k < n$) stuck on the floor with zero energy, i. e. $q_1 = \cdots = q_k = 0$, $p_1 = \cdots = p_k = 0$. Then the only natural way of defining the collisions of the orbit $\{\psi^t(x)|\, t \in \mathbb{R}\}$ is such that these particles will remain standing still forever on the floor with zero energy. We call these trajectories the degenerate ones.

**Proposition A.1.** *Suppose that the trajectory $\{\psi^t(x)|\, t \in \mathbb{R}\}$ is nondegenerate. Then the collision moments of this orbit can not accumulate at any real number $t$. In other words, there can only be finitely many collisions in finite time.*

**Proof.** Assume the opposite, i. e. that there is a (say, positive) accumulation point of collision moments. Denote by $t_0$ the smallest one of such positive accumulation points. Then all collision moments of the open interval $(0, t_0)$ form an increasing sequence $0 < t_1 < t_2 < \cdots < t_0$ such that $\lim_{m \to \infty} t_m = t_0$. Denote, as usual, by $\sigma_k = (i_k, i_k + 1)$ $(0 \leq i_k \leq n-1)$ the type of the collision taking place at time $t_k$, where $(0, 1)$ means the collision with the floor. We say that an integer $i$ $(0 \leq i \leq n-1)$ is *essential* if and only if there are infinitely many natural numbers $k$ with $i_k = i$. By classically known results ([G(1978)], [V(1979)], [G(1981)]), without the floor collision there can only be finitely many collisions in finite time. Therefore, the set of essential indexes $i$ has the form $\{0, 1, \ldots, a\}$, $0 \leq a \leq n - 1$. We can assume that the origin $t = 0$ is already chosen in such a way that $i_k \leq a$ for every positive integer $k$. We are now focusing on the limit $\lim_{t \nearrow t_0} v_i(t)$ of the velocities $v_i(t)$ of the particles as $t \nearrow t_0$.

**Lemma A.2.** *For every $i$, $1 \leq i \leq a + 1$, the limit*

$$\lim_{t \nearrow t_0} v_i(t) := v_i^-(t_0)$$

*exists.*

**Proof.** We start with the case $i = a + 1$. Since $\dot{v}_{a+1} = -1$ between collisions and $v_{a+1}^+(t_k) > v_{a+1}^-(t_k)$ for $i_k = a$, we conclude that

$$(A.3) \qquad \lim_{t \nearrow t_0} v_{a+1}(t) = v_{a+1}(0) - t_0 + \sum_{k=1}^{\infty} \left[ v_{a+1}^+(t_k) - v_{a+1}^-(t_k) \right].$$

This settles the case $i = a + 1$.

Suppose now that $1 \leq i \leq a$ and the lemma has been proved for $i+1, \ldots, a+1$. Then again we have $\dot{v}_i = -1$ between collisions of type $i-1$ or $i$ and $v_i^+(t_k) > v_i^-(t_k)$ if $i_k = i - 1$, while $v_i^+(t_k) < v_i^-(t_k)$ if $i_k = i$. Since the lemma is supposed to be valid for $i+1, \ldots, a+1$, one concludes that

$$\sum_{k,\, i_k = i} \left[ v_i^-(t_k) - v_i^+(t_k) \right] < \infty.$$

Thus, an argument similar to the one yielding (A.3) provides

$$(A.4) \quad \lim_{t \nearrow t_0} v_i(t) = v_i(0) - t_0 - \sum_{k,\, i_k = i} \left[ v_i^-(t_k) - v_i^+(t_k) \right] + \sum_{k,\, i_k = i-1} \left[ v_i^+(t_k) - v_i^-(t_k) \right].$$

Hence the lemma follows. □

Finally, since there are infinitely many collisions of type $(i, i+1)$ $(0 \leq i \leq a)$, we conclude that $v_i^-(t_0) \leq v_{i+1}^-(t_0)$ (following from the fact that $v_i^+(t_k) < v_{i+1}^+(t_k)$ for $i_k = i$) and, similarly, $v_{i+1}^-(t_0) \leq v_i^-(t_0)$. (Here we use the natural convention $v_0(t) = 0$.) Thus $v_i^-(t_0) = 0$ for $i = 1, \ldots, a+1$ and, therefore, $\lim_{t \nearrow t_0} q_i(t) = 0$, as well.

However, the subsystem $\{1, 2, \ldots, a+1\}$ can not just loose its positive energy! The obtained contradiction finishes the proof of Proposition A.1. □

**Remark.** As it has been shown in [B-F-K(1995)], the number of collisions in a unit time interval is bounded in any semi-dispersing billiard. Since a small perturbation of a degenerate trajectory clearly provides an arbitrarily high number of collisions in unit time, we conclude that there is no way to introduce any Riemannian metric in the configuration space $\mathbf{Q}$ in such a way that

(i) all sectional curvatures of that metric are nonpositive,

(ii) the smooth components of the boundary $\partial \mathbf{Q}$ are convex from inside $\mathbf{Q}$, and

(iii) the flow $\{\psi^t\}$ is induced by the geodesic flow in $\mathbf{Q}$ and by the usual law of reflections at the boundary $\partial \mathbf{Q}$.